\magnification 1200
\def\R{{\rm I\kern-0.2em R\kern0.2em \kern-0.2em}}
\def\N{{\rm I\kern-0.2em N\kern0.2em \kern-0.2em}}
\def\P{{\rm I\kern-0.2em P\kern0.2em \kern-0.2em}}
\def\B{{\rm I\kern-0.2em B\kern0.2em \kern-0.2em}}
\def\Z{{\rm I\kern-0.2em Z\kern0.2em \kern-0.2em}}
\def\C{{\bf \rm C}\kern-.4em {\vrule height1.4ex width.08em depth-.04ex}\;}

\def\D{{\Delta}}

\def\z{{\zeta}}
\def\cC{{\cal C}}

\def\cU{{\cal U}}

\def\degr{{\hbox{ deg }}}
\font\ninerm=cmr8
\ 
\vskip 6mm
\centerline {\bf DEGREE AND HOLOMORPHIC EXTENSIONS}
\vskip 4mm
\centerline{Josip Globevnik}
\vskip 4mm
{\noindent \ninerm ABSTRACT\ \ Let $D$ be a bounded convex domain 
in $\C^N,\ N\geq 2$. We prove that a continuous map $\Phi\colon\ bD\rightarrow\C^N$ 
extends holomorphically through $D$ if and only if for every polynomial map $P\colon\ 
\C^N\rightarrow\C^N$ such that $\Phi+P\not= 0$ on $bD$, the degree of $\Phi+P|bD$ is 
nonnegative. We also prove another such theorem for more general domains.
} 
\vskip 6mm
\bf 1.\ Introduction\rm
\vskip 2mm
Let $D\subset\R^n$ be a bounded open set and let $\Psi\colon\ bD\rightarrow 
R^n\setminus \{ 0\} $ be a continuous map.
Let $\tilde\Psi $ be a continuous extension of $\Psi $ to $\overline D$. Approximate $\tilde \Psi $ 
on $\overline D$ uniformly 
by a map $G$ smooth in a neighbourhood of $\overline D$ such that $G(bD)\subset \R^n\setminus \{ 0\}$. 
Perturbing $G$ slightly we may assume that the origin $0$ is a regular value of $G$ so $G^{-1}(0)\cap D$ 
is a finite subset of $D$ and each point in $G^{-1}(0)\cap D$ is a regular point of $G$. Let $\nu $ be the 
number of points in $G^{-1}(0)\cap D$ at which the derivative $DG$ preserves orientation minus the number of 
points in $G^{-1}(0)\cap D$ 
at which $DG$ reverses orientation. The number $\nu$ depends neither on the choice of the extension 
$\tilde\Psi $ of $\Psi $ nor on the choice of $G$ provided that $G$ approximates $\tilde \Psi $ on $\overline D$ 
well enough [D]. 
It is called the degree of $\Psi$, \ $\nu = \hbox{ deg }\Psi$. It is known that if $\{ \Psi _t ,\ 0\leq t\le 1\}$,
is a continuous family  
of continuous maps from $bD$ to $\R ^n\setminus\{ 0\}$ then $\degr \Psi_1 = \degr \Psi_0$. In the special case 
when $D\subset \C$ 
is a bounded domain with smooth boundary and $\Psi \colon \ bD\rightarrow \C\setminus \{ 0\}$ is a continuous 
function then 
$2\pi \degr \Psi $ equals the change of argument of $\Psi $ along $bD$. [D]

Let $D\subset \C^N$ be a bounded domain and suppose that $\Phi\colon\ bD\rightarrow \C^N\setminus\{ 0\} $ 
is a continuous map which extends 
holomorphically through $D$. Then $\degr \Phi \geq 0$. To see this, observe first that 
perturbing $\Phi $ slightly does not 
change the degree and implies that all the zeros of the holomorphic extension $\tilde\Phi $ 
of $\Phi $ are regular points of 
$\tilde\Phi $. Since $\tilde\Phi $ is holomorphic, at each regular 
point $a$ of $\tilde\Phi $ the derivative 
$(D\tilde\Phi )(a)$, a $\C$-linear map looked upon as a linear map from $\R^{2N}$ to  $\R^{2N}$,  
preserves orientation. In particular, $\degr \Phi $ is equal to the number of points $a\in D$ such 
that $\tilde \Phi (a)=0$ hence $\degr \Phi \geq 0$. 

Assume that $\Psi \colon bD\rightarrow \C^N$ is a continuous map. If $\Psi $ extends holomorphically through $D$ then 
by the preceding discussion $\degr (\Psi + F)\geq 0$ for every continuous map $F\colon \ bD\rightarrow \C^N$ 
that extends holomorphically through $D$ and is such that $\Psi +F\not=0$ on $bD$. It is 
known that the converse is true if $D$ is a smoothly bounded domain in $\C$.  
In the present paper we prove the converse for a large class of domains in $\C^N,\ N\geq 2$.

\vskip 3mm
\bf 2.\ The main results\rm
\vskip 2mm
Our main results are the following two theorems.  
\vskip 2mm
\noindent \bf THEOREM 2.1\ \it Let $D$ be a bounded convex domain in $\C^N,\ N\geq 2$. A 
continuous map $\Phi\colon\ bD\rightarrow \C^N$ extends holomorphically through $D$ if and 
only if for every polynomial map $P\colon\ \C^N\rightarrow\C^N$ such that $\Phi +P\not= 0 $
on $bD$, the degree of $\Phi +P|bD$ is nonnegative.
\vskip 2mm
\noindent \bf THEOREM 2.2\ \it Let $N\geq 2$ and let $D\subset \C^N$ be a bounded domain with 
$\cC^2$ boundary such that $\overline D$ has a Stein neighbourhood basis. A continuous 
map $\Phi\colon \ bD\rightarrow \C^N$ 
extends holomorphically through $D$ if and only if for each holomorphic map $G$ from a 
neighbourhood of $\overline D$ (that may depend on $G$) to $
\C^N$ such that $\Phi + G\not= 0$ on $bD$, the degree of $\Phi + G|bD$ is nonnegative. 
\vskip 2mm \rm
The theorems are known in the case when $N=1$: Given a bounded domain $D$ in $\C$ let $A(D)$ be 
the algebra of all continuous functions on $\overline D$ which 
are holomorphic on $D$. If $bD$ consists of finitely many pairwise disjoint 
simple closed curves then a continuous function $\Phi $ on 
$bD$ extends holomorphically through $D$ if and only if for each $G\in A(D)$ 
such that $\Phi + G\not= 0$ on $bD$, the change of argument of $\Phi +G$
along $bD$ is nonnegative [G2]. In fact, $A(D)$ here may be replaced by any dense subset 
of $A(D)$, for instance, by the set of functions holomorphic in a neighbourhood 
of $\overline D$ (which may depend on the function) [S2]. In particular, 
if $D$ is simply connected, it suffices to take for $G$ the polynomials.

\vskip 4mm
\bf 3.\ The degree of a special map on the intersection of $bD$ with a complex line \rm
\vskip 2mm
Denote by $L$ the $z_1$-axis in $\C^N,\ N\geq 2$,$$
L=\{ (\z,0,\cdots, 0)\colon \ \z\in\C\}.
$$
\vskip 1mm
\noindent\bf PROPOSITION 3.1\ \it Let $D\subset \C^N$ be a bounded domain. Suppose that $L$ 
meets $D$ and that $L\cap bD$ is the boundary of $L\cap D$ in $L$. Let 
$\Omega = \{ \z\in\C\colon \ (\z ,0,\cdots ,0)\in D\}$ so $b\Omega = \{ \z\in\C\colon \ (\z ,0,\cdots ,0)\in bD\}$.
Let $\varphi\colon\ bD\rightarrow \C$ be a continuous function such that $\varphi (\z,0,\cdots 0)\not= 0\ 
(\z \in b\Omega )$. 
Define a continuous map $\Phi \colon\ bD\rightarrow \C^N\setminus\{ 0\}$ by
$$
\Phi (z_1,\cdots, z_N)=\bigl(\varphi (z_1,\cdots ,z_N), z_2,\cdots ,z_N\bigr)\ \ \bigl((z_1,\cdots ,z_N)\in bD\bigr).
$$
Then $\degr \Phi $ equals the degree of the map $\z\mapsto \varphi (\z ,0,\cdots ,0)\ (\z\in b\Omega )$. \rm
\vskip 2mm 
\noindent \bf Proof.\ \rm We first show the following;

\noindent There is an $\epsilon >0$ such that whenever $\varphi _1 $ is a continuous function on $bD$ such that 
$|\varphi _1(z)-\varphi (z)|<\epsilon\ (z\in L\cap bD)$ then the degrees of the maps 
$\z\mapsto\varphi (\z ,0,\cdots ,0)\ (\z\in b\Omega )$
 and  $\z\mapsto\varphi _1 (\z ,0,\cdots ,0)\ (\z\in b\Omega )$ are the same and, moreover, if
 $\Phi _1(z) = (\varphi _1(z), z_2,\cdots , z_N)\ (z\in bD)$ then $\degr \Phi _1 = \degr \Phi$. \rm 
 
 \noindent To see this, recall first that by our assumption, \ $\varphi (z)\not= 0\ (z\in L\cap bD)$ so
 there is an 
 $\epsilon >0$ such that if $\varphi _1\colon\ bD\rightarrow \C$ is a continuous function such that 
 $|\varphi _1 -\varphi | <\epsilon $ on $L\cap bD$ then $(1-\lambda )\varphi +\lambda \varphi _1\not= 0$ 
 on $L\cap bD$ 
 for each $\lambda,\ 0\leq\lambda\leq 1$. In particular, \ $(1-\lambda)\varphi (\z , 0,\cdots ,0)+
 \lambda \varphi _1(\z ,0,\cdots ,0)\not= 0\ (\z\in b\Omega ,\ 0\leq\lambda\leq 1)$ which implies that the degrees of the maps 
$\z\mapsto\varphi (\z ,0,\cdots ,0)\ (\z\in b\Omega )$
 and  $\z\mapsto\varphi _1 (\z ,0,\cdots ,0)\ (\z\in b\Omega )$ are the same. Fix such $\varphi _1 $ and let 
 $\Phi _1 (z)= 
 (\varphi _1(z), z_2, \cdots, z_N)$. Consider
 $(1-\lambda )\Phi (z)+\lambda \Phi _1(z)=((1-\lambda )\varphi (z)+\lambda \varphi _1 (z),
 z_2, \cdots ,z_N)$. If $z\in L\cap bD$ then $(1-\lambda )\varphi (z)+\lambda \varphi _1(z)\not=0$ so 
 $(1-\lambda )\Phi (z)+\lambda \Phi _1(z)\not= 0\ (0\leq\lambda\leq 1)$. If $z\in bD\setminus L$ 
 then $(z_2, \cdots z_N)\not= 0$ so again $(1-\lambda )\Phi (z)+\lambda \Phi _1(z)\not= 0\ (0\leq\lambda\leq 1)$. 
 Thus, 
 $\Psi _\lambda = (1-\lambda )\Phi +\lambda \Phi _1,\ 0\leq\lambda\leq 1$, is a continuous family 
 of continuous maps
 from $bD$ to 
 $\C^N\setminus \{ 0\}$ so $\degr \Phi =\degr \Psi _0 =\degr \Psi_1 =\degr \Phi_1$. The statement is proved.
 
 Choose a smooth complex valued function $\omega $ on $\C $ which satisfies
 $$
 |\omega (\z )-\varphi (\z ,0, \cdots, 0)|<\varepsilon \ \ \ (\z\in b\Omega )
 $$
 and is such that $0$ is its regular value. Define a smooth function $\varphi _1$ on $\C^N$ by
 $$
 \varphi _1(z_1, z_2, \cdots ,z_N) = \omega (z_1)
 $$
 and define  $\Phi _1(z)=(\varphi _1(z), z_2,\cdots z_N)$. By the preceding paragraph the proof of
 Proposition 3.1 will be complete 
 once we have shown that $\degr \Phi _1 $ is the same as the degree of the 
 map $\z\mapsto \omega (\z )\ (\z\in b\Omega )$. 
 
 By the assumption, $\omega (x+iy)= u(x,y)+iw(x,y)= (u(x,y),w(x,y))$ has finitely many zeros 
 $a_j= p_j+iq_j, \ 1\leq j\leq m,$ in $\Omega $ and each of these zeros is a regular point of $\omega$. Moreover, 
 by the construction, the map $\Phi _1$ has precisely the zeros $(a_j,0,\cdots,0), \ 1\leq j\leq m$, in $D$.
 
 Suppose that $a=p+iq= (p,q)$ is one of the zeros of $\omega $  so $(a,0, \cdots 0)$ is 
 a zero of $\Phi _1$. Since $a$ is a regular point of $\omega $ the derivative 
 $$
 (D\omega )(a) =\left[
 \eqalign{{{\partial u}\over{\partial x}}(p,q), {{\partial u}\over{\partial y}}(p,q)\cr
 {{\partial w}\over{\partial x}}(p,q), {{\partial w}\over{\partial y}}(p,q)\cr}\right],
 $$
 looked upon as a linear map from $\R^2$ to $\R^2$, is 
 nonsingular. Since $\Phi _1(x_1,y_1,\cdots ,x_N, y_N)= (u(x_1,y_1), w(x_1,y_1), x_2, y_2,\cdots ,
 x_N, y_N) $ the 
 derivative of $\Phi _1$ at $(a,0,\cdots ,0)=$ $ (p,q,0,$ $\cdots ,$ $0)$, looked upon as a 
 linear map from $\R^{2N}$ 
 to $\R^{2N}$ is
 $$
 (D\Phi_1)(a,0,\cdots ,0)= \left[
 \eqalign{{{\partial u}\over{\partial x_1}}(p,q), {{\partial u}\over{\partial y_1}}(p,q),\ \ &O\cr
 {{\partial w}\over{\partial x_1}}(p,q), {{\partial w}\over{\partial y_1}}(p,q),\ \ &O \cr
 \ \ \ O, \ \ \ \ \ \ \  \ \ O,\ \  &I \cr}\right] .
 $$
where $I$ is the identity matrix of order $2N-2$. Thus, $(D\Phi _1)(a,0,\cdots ,0)$ is also nonsingular and 
$\hbox{det} (D\omega )(a) = \hbox{det}(D\Phi _1)(a,0,\cdots, 0)$. It follows that the maps $(D\omega )(a)$ and 
$ (D\Phi _1)(a,0,\cdots, 0)$ either both preserve orientation or both reverse orientation. Since
$\Phi_1(z)=0$ for $z\in D$ if and only if $z=(a_j, 0,\cdots, 0)$ for some $j,\ 1\leq j\leq m$, 
it follows that the degree 
of $\z \mapsto \omega (\z )\ (\z\in b\Omega )$ equals $\degr \Phi _1$. This completes the proof.
\vskip 1mm
\noindent We shall also need
\vskip 2mm
\noindent\bf PROPOSITION 3.2\ \it Let $\Phi = (\Phi _1,\cdots , \Phi _N)$ be a continuous map 
from $bD$ to $\C^N\setminus \{ 0\}$. 
Let $t_j>0\ (1\leq j\leq N)$ and let $\Psi = (t_1\Phi _1,\cdots ,t_N\Phi _N)$. Then 
$\degr \Psi = \degr \Phi $. \rm
\vskip 2mm
\noindent\bf Proof.  \rm $\Theta _\lambda = (1-\lambda )\Phi + \lambda \Psi, \ 0\leq \lambda\leq 1$ 
is a continuous family of maps from $bD$ to $\C^N\setminus\{ 0\}$ such 
that $\Theta _0=\Phi,\ \ \Theta _1=\Psi $.  It follows that $\degr \Psi =\degr \Phi $. This completes the proof. 
\vskip 4mm
\bf 4.\ Proofs of Theorems 2.1 and 2.2 \rm
\vskip 2mm
\noindent\bf Lemma 4.1\ \it 
Let $D\subset \C $ be a bounded open set with $\cC^1$ boundary. A continuous function $\Phi $ on $bD$ 
extends holomorphically through $D$ if 
and only if for each function $G$, holomorphic in a neighbourhood of $\overline D$, 
such that \ $\Phi + G\not= 0$ on $bD$, the degree of 
$\Phi+G|bD$ is nonnegative.  \rm
\vskip 1mm
\noindent \bf Proof.\ \rm Observe first that $D=D_1\cup\cdots\cup D_m$ 
where $D_j,\ 1\leq j\leq m$, 
are domains with pairwise disjoint closures and each $bD_j,\ 1\leq j\leq m$, consists of 
finitely many pairwise disjoint simple 
closed curves. The only if part follows from the argument principle. To prove the if part, 
assume that $\Phi\colon\ bD\rightarrow\C$ 
is a continous function that does not extend holomorphically through $D$. So for some $j$, 
the function 
$\Phi|bD_j$ does not extend holomorphically through $D_j$ which implies [G2] that there is a 
function $H\in A(D_j)$ such that $\Phi + H\not=0 $ on $bD_j$ and that $\degr (\Phi |bD_j+H|bD_j)$ 
is negative. Since $H$ 
can be approximated on $\overline{D_j}$ arbitrarily well by rational functions with poles outside 
$\overline{D_j}$ \ [S2] we may assume 
that $H$ is holomorphic on a neighbourhood $U$ of $\overline{D_j}$ whose closure misses 
$bD_k,\ 1\leq k\leq m,\ k\not= j$. Adding a sufficiently large constant 
$T_k$ to $H|bD_k,\ k\not= j, 1\leq k\leq N$,  will make
the degree 
of $H|bD_k+T_k$ equal zero. So putting $ H\equiv T_k$ on $bD_k,\ 1\leq k\leq m,\ k\not= j$, we get 
a function $H$, holomorphic on 
a neighbourhood of $\overline D$ such that $\Phi + H\not= 0$  on $bD$ and such that
$\degr (\Phi+H|bD)$ is negative.  This proves the if part and completes the 
proof.
\vskip 2mm
\noindent\bf Proof of Theorem 2.2.\ \rm The only if part was proved at the end of Section 2. To prove the if
part assume that
$\Phi = (\Phi_1, \cdots ,\Phi _N)$ does not extend holomorphically through $D$. Then one of the components, 
say $\Phi_1$, 
does not extend holomorphically through $D$ which implies that there is a complex line $L$ meeting $D$ and 
meeting $bD$ transversely such that 
$\Phi_1|(L\cap bD)$ does not extend holomorphically through $L\cap D$\ [GS]. After a translation and rotation we
may assume 
with no loss of generality that $L$ is the $z_1$-axis. Let $\Omega =\{\z \in\C\colon\ (\z, 0,\cdots ,0\}\in D\}$,
so  
$b\Omega = \{\z\in\C\colon\ (\z, 0, \cdots ,0)\in bD\}$. The function $\z\mapsto \Phi_1(\z , 0,\cdots ,0)$ is 
continuous on $b\Omega$ and does 
not extend holomorphically through $\Omega $ which, by Lemma 4.1 implies that there is a holomorphic function $g$ 
on an open 
neighbourhood $U$ of $\overline\Omega$ in $\C$ such that
$\Phi _1(\z, 0,\cdots ,0)+g(\z )\not= 0 \ (\z\in b\Omega )$ 
and such that the map
$\z\mapsto \Phi _1(\z, 0,\cdots ,0)+g(\z )\ (\z\in b\Omega)$ has negative degree.
Since $\overline D$ has 
a Stein neighbourhood basis it has arbitrarily small pseudoconvex neighbourhoods so there is a 
pseudoconvex domain $\Sigma $ containing $\overline D$ such that 
$\{(\z ,0,\cdots ,0)\colon\ \z\in U\}\cap\Sigma $ is a 
closed subset of $\Sigma$ and hence a closed one dimensional submanifold of
the pseudoconvex domain $\Sigma$. It follows
[GR, p.\ 245] that there is a holomorphic function $H_1$ on $\Sigma $ such that
$$
H_1(\z, 0,\cdots ,0)=g(\z )\ \ ((\z,0,\cdots 0)\in\Sigma).
$$
Thus, $H_1$ is holomorphic in a neighbourhood of $\overline D$. 

Write $z=(z_1,\cdots,z_N)$. By Proposition 3.1  the map from $bD$ to $\C^N\setminus \{ 0\}$ given by 
$$
z\mapsto (\Phi _1(z)+H_1(z), z_2, \cdots ,z_N) \ \ (z\in bD) 
\eqno (4.1)
$$
has the same degree as the map from $b\Omega $ to $\C\setminus\{ 0\}$ given by 
$$
\z\mapsto \Phi_1(\z ,0,\cdots, 0)+g(\z )\ \ (\z\in b\Omega )
$$
which implies that the degree of the map (4.1) is negative.

Perturbing the map (4.1) slightly will not change the degree so one can choose $T>0$ so large that the map 
$$
z\mapsto \bigl(\Phi_1(z)+H_1(z), z_2+ \Phi_2(z)/T,\cdots ,z_N+\Phi _N(z)/T\bigr)
\eqno (4.2)
$$
maps $bD$ into $\C^N\setminus \{ 0\}$ and has the same degree as the map (4.1). 
By Proposition 3.2 the degree of the map
$$
z\mapsto (\Phi_1(z)+H_1(z), \Phi _2(z) + Tz_2, \cdots ,\Phi _N(z) + Tz_N)\ \ (z\in bD)
$$
is the same as the degree of the map (4.2). So, setting  $H_j(z)\equiv Tz_j\ (z\in U, \ 2\leq j\leq N)$
we have constructed a holomorphic map $H\colon\ U\rightarrow \C^N$ such that $\Phi +H\not= 0$ on $bD$ and such that $\degr (\Phi + H|bD) $
is negative. This completes the proof of Theorem 2.2.
\vskip 2mm \noindent\bf 
Proof of Theorem 2.1\ \rm If $D$ is convex then $\Omega $ in the proof, being convex, is a simply connected 
domain so $g$ can be chosen to be a polynomial and for $H_1$ one can take a polynomial on $\C^N$ defined by 
$H_1(z_1,\cdots, z_N)=g(z_1)$ to have $H_1(\z ,0, \cdots ,0)=g(\z )\ \ ((\z, 0,\cdots, 0)\in\overline D)$.  
One finishes the proof as the proof of Theorem 1.2. 
Theorem 2.1 is proved. 
\vskip 4mm
\bf 5.\ Consequences and remarks \rm
\vskip 2mm
A continuous function $\Phi \colon\ bD\rightarrow\C$ extends holomorphically through 
$D$ if and only if the map $z\mapsto (\Phi (z),0,\cdots, 0)$ extends holomorphically through $D$ which gives
\vskip 2mm
\noindent\bf COROLLARY 5.1\ \it Let $N\geq 2 $ and let $D\subset \C^N$ be a bounded domain with $\cC^2$ boundary
such that $\overline D$ has a Stein neighbourhood basis. A continuous function $\Phi\colon\ bD\rightarrow C$ 
extends holomorphically through $D$ if and only if for any $N$-tuple of functions $G_1,\cdots ,G_N$ holomorphic 
in a neighbourhood of $\overline D$ and such that $(\Phi+G_1, G_2,\cdots, G_N)\not= 0$ on $bD$, the degree of the map 
$z\mapsto (\Phi (z)+G_1(z), G_2(z),\cdots ,G_N(z))\ (z\in bD)$ is nonnegative. 
\vskip 2mm
\noindent \rm This strenghtens the main result of [S1] where it was shown that $\Phi $ 
extends holomorphically through $D$ if and only if the degree 
of $z\mapsto (H_1(z,\Phi (z),\cdots ,H_N(z,\Phi (z)))\ (z\in bD)$ is nonnegative whenever $H_1,\cdots ,H_N$ are 
holomorphic functions on a neighbourhood of $\overline D\times\C$ in $\C^{N+1}$ such that
$(H_1(z,\Phi (z),\cdots ,H_N(z,\Phi (z)))\not= 0 \ (z\in bD)$.
\vskip 2mm
\noindent\bf COROLLARY 5.2\ \it Let $D\subset \C^N$ be a convex domain. A continuous function 
$\Phi\colon\ bD\rightarrow \C$ extends holomorphically through $D$ if and only if for every 
$N$-tuple of polynomials 
$P_1,\cdots ,P_N\colon\ \C^N\rightarrow \C$ such that $(\Phi+P_1, P_2, \cdots ,P_N)\not= 0$ on $bD$, 
the degree of the map 
$z\mapsto (\Phi (z)+P_1(z), P_2(z),\cdots ,P_N(z))\ (z\in bD)$ is nonnegative.
\vskip 2mm\rm
Theorem 2.1 and Corollary 5.2 hold for more general domains:

\noindent Suppose that $\cU$ is an open subset of the set of all complex lines in $\C^N,\ N\geq 2$, 
passing through the origin, and let 
$D\subset \C^N$ be a bounded domain with $\cC^2$ boundary such that $L\cap bD$ is connected for every complex
line $L$ of the form $L=z+\Sigma $, $z\in\C^n,\ \Sigma\in\cU$, which intersects $bD$ transversely. Then the 
statement of Theorem 1.1 holds.  To see this, notice first that the assumptions imply that $bD$ is 
connected. Further, if 
$\Phi $ is a continuous function on $bD$ such that for each $L$ as above, $\Phi | L\cap bD$ extends 
holomorphically through 
$L\cap D$ then $\Phi $ satisfies weak tangential Cauchy Riemann equations on $bD$ [GS] and hence, since $bD$ 
is connected, 
$\Phi $ extends holomorphically through $D$ [K]. Thus, if $\Phi $ is a continuous function on $bD$ that does not
extend 
holomorphically through $D$ then there is a complex line $L$ as above such that $\Phi| L\cap bD$ does not 
extend holomorphically through $L\cap D$.
Proceeding as in the proof of Theorem 2.2 we notice that the connectedness of $L\cap bD$ implies 
that $\Omega $ is a simply connected domain and we 
can finish the proof as the proof of Theorem 2.1. 

If $A\colon\ \C^N\rightarrow \C^N$ is an invertible $\R$-linear map and $D$ is a bounded domain that
contains the origin 
then $\degr (A|bD)$ equals $+1$ if $A$ preserves orientation on $\R^{2N}=\C^N$ and $-1$ if  $A$ reverses 
orientation on $\R^{2N}=\C^N$ . Recall that 
an invertible $\C$-linear map preserves orientation. So, if $A\colon\ \C^N\rightarrow \C^N$ is a $\C$-linear 
map then $A+B$ preserves orientation whenever $B\colon\ \C^N\rightarrow \C^N$ is a $\C$-linear map such that $A+B$ 
is invertible. This property characterizes $\C$-linearity:
\vskip 2mm
\noindent\bf PROPOSITION 5.3\ \it Let $A\colon\ \C^N\rightarrow \C^N$ be a $\R$-linear map such that $A+B$ 
preserves orientation whenever $B\colon\ \C^N\rightarrow\C^N$ is a $\C$-linear map such that $A+B$ is invertible. Then $A$ is 
$\C$-linear.
\vskip 1mm
\noindent\bf Proof.\ \rm Assume that $A\colon\ \C^N\rightarrow \C^N$ is
a $\R$-linear map which is not $\C$-linear. Write $A=(A_1,\cdots ,A_N)$ where $A_j,\ 1\leq j\leq N$, are 
$\R$-linear functionals at least one of which, say $A_1$, is not $\C$-linear. It follows that there is a complex line $L$ 
through the origin such that $A_1|L$ is not $\C$-linear. With no loss of generality assume that $L$ is the $z_1$-axis. Thus, 
$A_1((\z,0,\cdots ,0)) = \alpha \z + \beta \overline\z\ (\z\in\C)$ where $\beta\not= 0$. 
There is $T>0$ so large that 
the map from $\C^N$ to $\C^N$ given by
$$
z\mapsto \bigl(\beta\overline{z_1}, Tz_2+tA_2(z),\cdots ,Tz_N+tA_N(z)\bigr)
\eqno (5.1)
$$
is invertible for each $t,\ 0\leq t\leq 1$. Since the map (5.1), looked upon as a linear map 
from $\R^{2N}$ to $\R^{2N}$ is invertible for each $t,\ 0\leq t\leq 1$, and 
reverses orientation for $t=0$ it
follows that it reverses orientation for each 
$t, \ 0\leq t\leq 1$. In particular, if we define
$$
H(z) = (-\alpha z_1, Tz_2,\cdots , Tz_N)q \ (z\in\C^N)
$$
it follows that
$$
z\mapsto (A+H)(z)=(\beta \overline {z_1}, A_2(z)+Tz_2,\cdots , A_N(z)+Tz_N)
$$
is an invertible map which reverses orinetation. This completes the proof.
\vskip 2mm
Proposition 5.3 can be viewed as the simplest case of Theorem 2.1. It shows that for a small class 
of maps $\Phi $ - $\R$-linear maps - a small class of holomorphic maps $P$ - $\C$-linear maps - 
is needed to check 
the holomorphic extendibility. It is an obvious question whether one can go further and ask whether for 
the set of all polynomial maps in $z_1,\cdots ,z_N, \overline{z_1},\cdots \overline {z_N}$ of degree $\leq m$, 
to check the holomorphic extendibility through $D$ it is enough to take for $P$ the holomorphic polynomials 
of degree $\leq m$. We prove that this is the case when $D$ is a ball:
\vskip 2mm\noindent\bf PROPOSITION 5.4\ \it Let $D\subset \C^N$ be an open ball, let $m\in\N$ and let 
$\Phi\colon\ \C^N\rightarrow \C^N$ be a polynomial map in 
$z_1,\cdots ,z_N,\overline{z_1},\cdots, \overline{z_N}$ 
of degree $\leq m$. If $\degr (\Phi + P)|bD$ is nonnegative whenever $P\colon\ \C^N\rightarrow \C^N$ is a holomorphic polynomial of 
degree $\leq m$ such that $\Phi + P\not=0$ on $bD$, then $\Phi|bD$ extends holomorphically through $D$ (as a 
holomorphic polynomial of degree $\leq m$). 
\vskip 1mm
\noindent\bf Proof.\ \rm Denote $\D (a,r)=\{\z\in\C\colon \ |\z-a|<r\}$. Note first that the fact that
$\Phi\colon\ \C^N\rightarrow \C^N$ is a polynomial map in $z_1,\cdots ,z_N,
\overline{z_1},\cdots ,\overline{z_N}$ of degree $\leq m$ is invariant with respect to affine $\C$-linear 
change of coordinates. Note also that if $p\colon\ \C\rightarrow\C$ is a polynomial of degree $\leq m$ in 
$\z $ and $\overline\z$ 
then given $a\in\C$ and $r>0$ there are polynomials $p$ and $s$ of degree $\leq m$ such that
$$
p(\z ,\overline\z )= q(\z-a)+\overline{s(\z -a)}\ \ (\z \in b\D (a,r)). 
$$
Assume that $\Phi\colon\ \C^N\rightarrow\C^N$ is a polynomial of degree $\leq m$ in $z_1,\cdots ,z_N,
\overline{z_1},\cdots ,\overline{z_N}$. Then for every $Z,W\in\C^N$, $\z\mapsto \Phi (Z+\z W)\colon\ \C\rightarrow \C^N$ is a 
polynomial of degree $\leq m$  in $\z,\overline\z $. 

Suppose that $\Phi |bD$ does not extend holomorphically through $D$. Then one of the 
components $\Phi_1,\cdots, \Phi_N$, say 
$\Phi_1$, does not extend from $L\cap bD$ holomorphically through  $L\cap D$. After a translation and a 
unitary change of coordinates we may assume that $L$ is the $z_1$-axis. Then $L\cap D$ is an open disc so 
$\{\z\in\C\colon\ (\z, 0,\cdots, 0)\in D\}$ is an open disc, $\Omega = \D (a,r)$ and $b\Omega = b\D (a, r) $ 
is 
a circle. Now, $\z\mapsto \Phi _1(\z ,0,\cdots ,0)$ is a complex valued polynomial in $\z $ and 
$\overline\z$ of degree $\leq m$ so
$$
\Phi_1 (\z ,0,\cdots ,0)= q(\z-a) + \overline{s(\z -a)}\ \ (\z\in b\Omega)
$$
where $q$ and $s$ are špolynomials of degree $\leq m$. Since $\Phi_1 |L\cap bD$ does not extend holomorphically 
through $L\cap D$ it follows that the polynomial $s$ is nonconstant so there is a $b\in\C$ such that the 
change of argument of $\z\mapsto\overline{s(\z -a)}-b$ along $b\Omega $ is negative. It follows that the degree of
$$
\z\mapsto \Phi_1(\z ,0,\cdots ,0)-q(\z -a) -b\ \ (\z\in b\Omega)
$$
is negative. Define 
$$
P_1 (z_1,\cdots, z_N) = -q(z_1 -a)-b.
$$
Then $P_1\colon\ \C^N\rightarrow \C$ is a holomorphic polynomial of degree $\leq m$ and the degree of the map 
$$
\z\mapsto\Phi_1(\z ,0,\cdots, 0)+P_1(\z ,0,\cdots , 0)\ \ (\z\in b\Omega)
\eqno (5.2)
$$
is negative. As in the proof of Theorem 2.2 we see that there is $T>0$ so large that if $P_j(z)=Tz_j\ (2\leq j\leq N)$ and 
if $P=(P_1,\cdots, P_N)$ then $\Phi +P\not= 0$ on $bD$ and the degree of the map $z\mapsto \Phi (z)+P(z)\ (z\in bD)$ 
coincides with the degree of the map (5.2). Thus we have constructed a holomorphic polynomial map 
$P\colon\ \C^N\rightarrow\C^N$ of degree $\leq m$ such that $\Phi+P\not=0$ on $bD$ and such that 
the degree of $(\Phi +P)|bD$ is negative. This completes the proof.
\vskip 5mm
This work was supported 
in part by the Ministry of Higher Education, Science and Technology of Slovenia  
through the research program Analysis and Geometry, Contract No.\ P1-0291. 
\vfill
\eject
\centerline{\bf REFERENCES}
\vskip 5mm
\vskip 5mm

\noindent [AW]\ H.\ Alexander and J.\ Wermer: Linking 
numbers and boundaries of varieties. 

\noindent Ann.\ Math.\ 151 (2000) 125-150
\vskip 2mm
\noindent [D]\ K.\ Deimling:\ \it Nonlinear functional analysis.\rm \ Springer Verlag, Berlin,  1980
\vskip 2mm
\noindent [G1]\ J.\ Globevnik:\ Holomorphic extendibility and the argument principle.

\noindent Complex Analysis and Dynamical Systems II. Contemp.\ Math.\ 382 (2005) 171-175
\vskip 2mm
\noindent [G2]\ J.\ Globevnik:\ The argument principle and holomorphic extendibility.

\noindent Journ.\ d'Analyse.\ Math.\ 94 (2004) 385-395
\vskip 2mm
\noindent [GS]\ J.\ Globevnik, E.\ L.\ Stout:\ Boundary Morera theorems for holomorphic functions of several complex variables.

\noindent Duke Math.\ J.\ 64 (1991) 571-615
\vskip 2mm
\noindent [GP] V.\ Guillemin, A.\ Pollack:\ \it Differential topology.\rm 

\noindent Prentice-Hall, Englewood Cliffs, New Jersey 1974
\vskip 2mm
\noindent [GR]\ R.\ Gunning, H.\ Rossi:\ \it Analytic Functions of Several Complex Variables.\rm 

\noindent Prentice-Hall, Englewood Cliffs, New Jersey 1965
\vskip 2mm
\noindent [K]\ A.\ M.\ Kytmanov:\ \it The Bochner-Martinelli integral and its applications.\rm 

\noindent Birkhauser Verlag, Basel-Boston-Berlin 1995
\vskip 2mm
\noindent [R]\ R.\ M. Range:\ \it Holomorphic functions and integral 
representations in several complex variables. \rm

\noindent Springer-Verlag, New York-Berlin-Heidelberg-Tokyo, 1986
\vskip 2mm
\noindent [S1]\ E.L.Stout: Boundary values and mapping degree.

\noindent Michig.\ Math.\ J.\ 47 (2000) 353-368
\vskip 2mm
\noindent [S2]\ E.\ L.\ Stout: \it The Theory of Uniform Algebras.\ 

\noindent \rm Bogden and Quigley, 
Tarrytown -on-Hudson, N.Y. 1971
\vskip 2mm
\noindent [W] J.\ Wermer: The argument principle and boundaries of analytic varieties.

\noindent Oper. Theory Adv. Appl., 127, Birkhauser, Basel, 2001, 639-659

\vskip 10mm
\noindent Institute of Mathematics, Physics and Mechanics

\noindent University of Ljubljana, Ljubljana, Slovenia

\noindent josip.globevnik@fmf.uni-lj.si

\bye

\bye